\documentclass{article}

\usepackage{amssymb}
\newcommand{\res}{{\upharpoonright}}
\newcommand{\NN}{\mathbb{N}}
\newcommand{\cat}{{^\smallfrown}}
\newcommand{\limp}{\Rightarrow}
\newcommand{\leT}{\leq_{\mathrm{T}}}
\newcommand{\nleT}{\nleq_{\mathrm{T}}}
\newcommand{\eqT}{\equiv_{\mathrm{T}}}

\usepackage{amsthm}
\theoremstyle{definition}
\newtheorem{thm}{Theorem}[section]
\newtheorem{lem}[thm]{Lemma}
\newtheorem*{dfn}{Definition}
\newtheorem*{dfns}{Definitions}
\newtheorem*{ntn}{Notation}
\newtheorem{rem}{Remark}

\usepackage[backref=page,colorlinks=true,allcolors=blue]{hyperref}

\begin{document}

\title{Harrington's results on arithmetical singletons}

\author{Stephen G. Simpson\\
  Department of Mathematics\\
  Pennsylvania State University\\
  University Park, PA 16802, USA\\
  \href{http://www.math.psu.edu/simpson}{http://www.math.psu.edu/simpson}\\
  \href{mailto:simpson@math.psu.edu}{simpson@math.psu.edu}}

\date{First draft: November 8, 2012\\
  This draft: December 12, 2012}

\maketitle


\begin{abstract}
  We exposit two previously unpublished theorems of Leo Harrington.
  The first theorem says that there exist arithmetical singletons
  which are arithmetically incomparable.  The second theorem says that
  there exists a ranked point which is not an arithmetical singleton.
  Unlike Harrington's proofs of these theorems, our proofs do not use
  the finite- or infinite-injury priority method.  Instead they use an
  oracle construction adapted from the standard proof of the Friedberg
  Jump Theorem.
\end{abstract}

%
%
%
 


\section{Introduction}
\label{sec:intro}

\begin{dfns}
  Let $\NN=\{0,1,2,\ldots,n,\ldots\}=\{$the natural numbers$\}$.  We
  work in the Baire space $\NN^\NN$.  Points in $\NN^\NN$ are denoted
  $X,Y,Z,\ldots$ and sets in $\NN^\NN$ are denoted $P,Q,\ldots$.  A
  point $X$ or a set $P$ is said to be \emph{arithmetical} if it is
  $\Pi^0_n$ for some $n$, and \emph{arithmetical relative to $Y$} if
  it is $\Pi^{0,Y}_n$ for some $n$.  See for instance Rogers
  \cite[Chapters 14--16]{rogers}.  Two points $X$ and $Y$ are said to
  be \emph{arithmetically equivalent} if each is arithmetical relative
  to the other, and \emph{arithmetically incomparable} if neither is
  arithmetical relative to the other.  An \emph{arithmetical
    singleton} is a point $X$ such that the singleton set $\{X\}$ is
  arithmetical.  A \emph{ranked point} is a point $X$ such that $X\in
  P$ for some countable $\Pi^0_1$ set $P$.
\end{dfns}

\begin{rem}
  It is well known that each arithmetical singleton is arithmetical
  relative to $0^{(\alpha)}$ for some recursive ordinal $\alpha$, and
  each such $0^{(\alpha)}$ is itself an arithmetical singleton.  See
  for instance Sacks \cite[Chapter II]{sacks-hrt}.
\end{rem}

\begin{rem}
  \label{rem:tanaka}
  Tanaka \cite{tanaka-arith} observed that for any arithmetical set
  $P$ we can find a $\Pi^0_1$ set $Q$ and a one-to-one correspondence
  $F:P\cong Q$ such that each $X\in P$ is uniformly arithmetically
  equivalent to $F(X)$.  It follows that every arithmetical singleton
  is arithmetically equivalent to a $\Pi^0_1$ singleton, every member
  of a countable arithmetical set is arithmetically equivalent to a
  ranked point, and every nonempty countable arithmetical set contains
  an arithmetical singleton.
\end{rem}

\begin{rem}
  The purpose of this paper is to exposit two previously unpublished
  theorems due to Harrington
  \cite{harrington-arith,harrington-mclaughlin} concerning
  arithmetical singletons.
  \begin{enumerate}
  \item There exist arithmetically incomparable arithmetical
    singletons.  Equivalently, there exist arithmetically incomparable
    $\Pi^0_1$ singletons.  See Theorem \ref{thm:incomp} below.
  \item There exists a point which belongs to a countable arithmetical
    set but is not an arithmetical singleton.  Equivalently, there
    exists a ranked point which is not an arithmetical singleton.  See
    Theorem \ref{thm:mclaughlin} below.
  \end{enumerate}
\end{rem}

\begin{rem}
  Harrington's theorems on arithmetical singletons may be reformulated
  so as to yield significant insights concerning definability over the
  natural number system $\NN,+,\times,=$.  Note first that $X$ is
  arithmetical if and only if $X$ is \emph{explicitly definable over
    $\NN,+,\times,=$}, and $X$ is an arithmetical singleton if and
  only if $X$ is \emph{implicitly definable over $\NN,+,\times,=$}.
  Prior to Harrington, two well known results concerning definability
  over $\NN,+\,\times,=$ were as follows.
  \begin{enumerate}
  \item There exists an $X$ which is implicitly definable over
    $\NN,+,\times,=$ but not explicitly definable over
    $\NN,+,\times,=$.  (Namely, let $X=0^{(\omega)}=$ the Tarski truth
    set for $\NN,+,\times,=$.  See Rogers \cite[Theorems 14-X and
    15-XII]{rogers}.)
  \item There exist $X$ and $Y$ such that $X\oplus Y$ is implicitly
    definable over $\NN,+,\times,=$ but neither $X$ nor $Y$ is
    implicitly definable over $\NN,+,\times,=$.  (Namely, let $X$ and
    $Y$ be Cohen generic over $\NN,+,\times,=$ such that $X\oplus
    Y\eqT0^{(\omega)}$.  See Feferman \cite{feferman-generic} or
    Rogers \cite[Exercise 16-72]{rogers}.)
  \end{enumerate}
  Note also that $X$ is arithmetical relative to $Y$ if and only if
  $X$ is \emph{explicitly definable over $\NN,+,\times,=,Y$}.  We then
  see that the following result due to Harrington is complementary to
  results 1 and 2.
  \begin{enumerate}
    \setcounter{enumi}{2}
  \item There exist $X$ and $Y$ such that
    \begin{enumerate}
    \item $X$ is implicitly definable over $\NN,+,\times,=$,
    \item $Y$ is implicitly definable over $\NN,+,\times,=$,
    \item $X$ is not explicitly definable over $\NN,+,\times,=,Y$,
    \item $Y$ is not explicitly definable over $\NN,+,\times,=,X$.
    \end{enumerate}
    (Namely, let $X$ and $Y$ be as in Theorem \ref{thm:incomp} below.)
  \end{enumerate}
\end{rem}

\begin{rem}
  Harrington \cite{harrington-arith,harrington-mclaughlin} and Gerdes
  \cite{gerdes-nusm} have applied the method of
  \cite{harrington-arith,harrington-mclaughlin} to prove many other
  remarkable theorems.  See for instance Remark \ref{rem:alpha} below.
  However, we choose not to present those applications here.  Instead
  we content ourselves with providing an accessible introduction to
  the method, including detailed proofs of two of the more striking
  theorems.
\end{rem}

\begin{rem}
  The plan of this paper is as follows.  In \S\ref{sec:simplified} we
  warm up by proving simplified versions of Harrington's theorems on
  arithmetical singletons.  In \S\ref{sec:full} we prove the full
  versions.
\end{rem}

\section{The simplified versions}
\label{sec:simplified}

\begin{dfns}
  Points $A,B\in\NN^\NN$ may be viewed as \emph{Turing oracles}.  We
  write $\{e\}^A(i)=j$ to mean that the $e$th Turing machine with
  oracle $A$ and input $i$ halts with output $j$.  We write
  $\{e\}^A(i)\downarrow$ (respectively $\uparrow$) to mean that the
  $e$th Turing machine with oracle $A$ and input $i$ halts (does not
  halt).  We write $A\leT B$ to mean that $A$ is \emph{Turing
    reducible to} $B$, i.e., $\exists e\,\forall
  i\,(A(i)=\{e\}^B(i))$.  We write $A\eqT B$ to mean that $A$ is
  \emph{Turing equivalent to} $B$, i.e., $A\leT B$ and $B\leT A$.  We
  define $A\oplus B\in\NN^\NN$ by the equations $(A\oplus B)(2i)=A(i)$
  and $(A\oplus B)(2i+1)=B(i)$.  We write $A'=$ the \emph{Turing jump}
  of $A$, defined by
  \[
  A'(e)=\left\{
    \begin{array}{ll}
      1 &\mbox{if }\{e\}^A(e)\downarrow\,,\\[4pt]
      0 &\mbox{if }\{e\}^A(e)\uparrow.
    \end{array}
  \right.
  \]
  We write $A^{(n)}=$ the $n$th Turing jump of $A$, defined
  inductively by letting $A^{(0)}=A$ and $A^{(n+1)}=(A^{(n)})'$.
  Recall that $A$ is arithmetical relative to $B$ if and only if
  $A\leT B^{(n)}$ for some $n$.
\end{dfns}

\begin{lem}
  \label{lem:0}
  Given a $\Pi^{0,A'}_1$ set $P$ we can find a $\Pi^{0,A}_1$ set $Q$
  and a homeomorphism $F:P\cong Q$ such that $X\oplus A\eqT F(X)\oplus
  A$ uniformly for all $X\in P$.
\end{lem}

\begin{proof}
  Since $P$ is a $\Pi^{0,A'}_1$ set, it follows that $P$ is a
  $\Pi^{0,A}_2$ set, say $P=\{X\mid\forall i\,\exists j\,R(X,i,,j)\}$
  where $R$ is an $A$-recursive predicate.  Define $F:P\cong Q=F(P)$
  by letting $F(X)=X\oplus\widehat{X}$ where $\widehat{X}(i)=$ the
  least $j$ such that $R(X,i,j)$ holds.  Clearly $Q$ is a
  $\Pi^{0,A}_1$ set and $X\oplus A\eqT F(X)\oplus A$ uniformly for all
  $X\in P$.
\end{proof}

\begin{lem}
  \label{lem:PQ}
  Given a $\Pi^{0,A'}_1$ set $P$ we can find a $\Pi^{0,A}_1$ set $Q$
  and a homeomorphism $H:P\cong Q$ such that $X\oplus A'\eqT
  H(X)\oplus A'\eqT(H(X)\oplus A)'$ uniformly for all $X\in P$.
\end{lem}

In order to prove Lemma \ref{lem:PQ}, we first present some general
remarks concerning strings, trees, and treemaps.

\begin{ntn}[strings]
  Let $\NN^*=\bigcup_{l\in\NN}\NN^l=$ the set of \emph{strings}, i.e.,
  finite sequences of natural numbers.  For $\sigma=\langle
  n_0,n_1,\ldots,n_{l-1}\rangle\in\NN^*$ we write $\sigma(i)=n_i$ for
  all $i<|\sigma|=l=$ the \emph{length} of $\sigma$.  For
  $\sigma,\tau\in\NN^*$ we write $\sigma\cat\tau=$ the
  \emph{concatenation}, $\sigma$ followed by $\tau$, defined by the
  conditions $|\sigma\cat\tau|=|\sigma|+|\tau|$,
  $(\sigma\cat\tau)(i)=\sigma(i)$ for all $i<|\sigma|$, and
  $(\sigma\cat\tau)(|\sigma|+i)=\tau(i)$ for all $i<|\tau|$.  We write
  $\sigma\subseteq\tau$ if $\sigma\cat\rho=\tau$ for some $\rho$.  If
  $|\sigma|\ge n$ we write $\sigma\res n=$ the unique
  $\rho\subseteq\sigma$ such that $|\rho|=n$.  If $|\sigma|=|\tau|=n$
  we define $\sigma\oplus\tau\in\NN^*$ by the conditions
  $|\sigma\oplus\tau|=2n$ and $(\sigma\oplus\tau)(2i)=\sigma(i)$ and
  $(\sigma\oplus\tau)(2i+1)=\tau(i)$ for all $i<n$.
\end{ntn}

\begin{dfn}[trees]
  A \emph{tree} is a set $T\subseteq\NN^*$ such that
  \begin{center}
    $\forall\rho\,\forall\sigma\,((\rho\subseteq\sigma$ and $\sigma\in
    T)\limp\rho\in T)$.
  \end{center}
  For any tree $T$ we write
  \begin{center}
    $[T]=\{$paths through $T\}=\{X\mid\forall n\,(X\res n\in T)\}$.
  \end{center}
\end{dfn}

\begin{rem}
  \label{rem:Pi01A}
  It is well known that the following statements are pairwise
  equivalent.
  \begin{enumerate}
  \item $P$ is a $\Pi^{0,A}_1$ set.
  \item $P=[T]$ for some $\Pi^{0,A}_1$ tree $T$.
  \item $P=[T]$ for some $A$-recursive tree $T$.
  \item $P=\{X\mid X\oplus A\in [T]\}$ for some recursive tree $T$.
  \end{enumerate}
\end{rem}

\begin{dfn}[treemaps]
  Let $T$ be a tree.  A \emph{treemap} is a function $F:T\to\NN^*$
  such that
  \begin{center}
    $ F(\sigma\cat\langle i\rangle)\supseteq F(\sigma)\cat\langle
    i\rangle$
  \end{center}
  for all $\sigma\in T$ and all $i\in\NN$ such that $\sigma\cat\langle
  i\rangle\in T$.  We then have another tree
  \begin{center}
    $F(T)=\{\tau\mid\exists\sigma\,(\sigma\in T$ and $\tau\subseteq
    F(\sigma))\}$.
  \end{center}
  Thus $P=[T]$ and $F(P)=[F(T)]$ are closed sets in the Baire space
  and we have a homeomorphism $F:P\cong F(P)$ defined by
  $F(X)=\bigcup_{n\in\NN}F(X\res n)$ for all $X\in P$.  Note also that
  the composition of two treemaps is a treemap.  A treemap
  $F:T\to\NN^*$ is said to be \emph{$A$-recursive} if it is the
  restriction to $T$ of a partial $A$-recursive function.
\end{dfn}

\begin{rem}
  \label{rem:FT}
  Let $T$ be a tree and let $F:T\to\NN^*$ be a treemap.  Given
  $\tau\in F(T)$ let $\sigma\in T$ be minimal such that $\tau\subseteq
  F(\sigma)$.  Then $\sigma$ is a \emph{substring} of $\tau$, i.e.,
  $\sigma=\langle\tau(j_0),\tau(j_1),\ldots,\tau(j_{l-1})\rangle$ for
  some $j_0<j_1<\cdots<j_{l-1}<|\tau|$.  Thus, in the definition of
  $F(T)$, the quantifier $\exists\sigma$ may be replaced by a bounded
  quantifier,
  \begin{center}
    $F(T)=\{\tau\mid(\exists\sigma$ substring of $\tau)\,(\sigma\in T$
    and $\tau\subseteq F(\sigma))\}$.
  \end{center}
  This implies that, for instance, if $F$ and $T$ are $A$-recursive
  then so is $F(T)$.
\end{rem}

We are now ready to prove Lemma \ref{lem:PQ}.

\begin{proof}[Proof of Lemma \ref{lem:PQ}]
  Given $A$ we construct a particular $A'$-recursive treemap
  $G:\NN^*\to\NN^*$.  We define $G(\sigma)$ by induction on $|\sigma|$
  beginning with $G(\langle\rangle)=\langle\rangle$.  If $G(\sigma)$
  has been defined, let $e=|\sigma|$ and for each $i$ let
  $G(\sigma\cat\langle i\rangle)=$ the least $\tau\supseteq
  G(\sigma)\cat\langle i\rangle$ such that $\{e\}_{|\tau|}^{\tau\oplus
    A}(e)\downarrow$ if such a $\tau$ exists, otherwise
  $G(\sigma\cat\langle i\rangle)=G(\sigma)\cat\langle i\rangle$.
  Clearly $G$ is an $A'$-recursive treemap, and our construction of
  $G$ implies that for all $e$ and $X$, $\{e\}^{G(X)\oplus
    A}(e)\downarrow$ if and only if $\{e\}_{|G(X\res e+1)|}^{G(X\res
    e+1)\oplus A}(e)\downarrow$.  Thus $X\oplus A'\eqT G(X)\oplus
  A'\eqT(G(X)\oplus A)'$ uniformly for all $X$.

  Let $G$ be the $A'$-recursive treemap which was constructed above.
  Let $P$ be a $\Pi^{0,A'}_1$ set.  By Remark \ref{rem:FT} we know
  that the restriction of $G$ to $P$ maps $P$ homeomorphically onto
  another $\Pi^{0,A'}_1$ set $G(P)$.  Applying Lemma \ref{lem:0} to
  $G(P)$ we obtain a $\Pi^{0,A}_1$ set $Q$ and a homeomorphism
  $F:G(P)\cong Q$ such that $Y\oplus A\eqT F(Y)\oplus A$ uniformly for
  all $Y\in G(P)$.  Thus $H=F\circ G$ is a homeomorphism of $P$ onto
  $Q$, and for all $X\in P$ we have $G(X)\oplus A\eqT F(G(X))\oplus
  A=H(X)\oplus A$ uniformly, hence $X\oplus A'\eqT H(X)\oplus
  A'\eqT(H(X)\oplus A)'$ uniformly, Q.E.D.
\end{proof}

\begin{rem}
  Our proof of Lemma \ref{lem:PQ} via treemaps is similar to the proof
  of \cite[Lemma 5.1]{masshyp}.  Within our proof of Lemma
  \ref{lem:PQ}, the construction of the specific treemap $G$ is the
  same as the standard proof of the Friedberg Jump Theorem as
  presented for instance in Rogers \cite[\S13.3]{rogers}.
\end{rem}

Let $0$ denote the constant zero function, so that $0^{(n)}=$ the
$n$th jump of $0$.

\begin{lem}
  \label{lem:Pn}
  Given a $\Pi^{0,0^{(n)}}_1$ set $P_n$ we can find a $\Pi^0_1$ set
  $P_0$ and a homeomorphism $H^n_0:P_n\cong P_0$ such that $X_n\oplus
  0^{(n)}\eqT X_0\oplus 0^{(n)}\eqT X_0^{(n)}$ uniformly for all
  $X_n\in P_n$ and $X_0=H^n_0(X_n)\in P_0$.
\end{lem}

\begin{proof}
  The proof is by induction on $n$.  For $n=0$ there is nothing to
  prove.  For the inductive step, given a $\Pi^{0,0^{(n+1)}}_1$ set
  $P_{n+1}$ apply Lemma \ref{lem:PQ} with $A=0^{(n)}$ to obtain a
  $\Pi^{0,0^{(n)}}_1$ set $P_n$ and a homeomorphism $H_n:P_{n+1}\cong
  P_n$ such that $X_{n+1}\oplus 0^{(n+1)}\eqT H_n(X_{n+1})\oplus
  0^{(n+1)}\eqT(H_n(X_{n+1})\oplus 0^{(n)})'$ uniformly for all
  $X_{n+1}\in P_{n+1}$.  Then apply the inductive hypothesis to $P_n$
  to find a $\Pi^0_1$ set $P_0$ and a homeomorphism $H^n_0:P_n\cong
  P_0$ such that $X_n\oplus 0^{(n)}\eqT X_0\oplus 0^{(n)}\eqT
  X_0^{(n)}$ uniformly for all $X_n\in P_n$.  Letting
  $H^{n+1}_0=H_n\circ H^n_0$ it follows that $X_{n+1}\oplus
  0^{(n+1)}\eqT X_0\oplus 0^{(n+1)}\eqT X_0^{(n+1)}$ uniformly for all
  $X_{n+1}\in P_{n+1}$ and $X_0=H^{n+1}_0(X_{n+1})\in P_0$, Q.E.D.
\end{proof}

We now use Lemma \ref{lem:Pn} to prove simplified versions of
Harrington's theorems.

\begin{thm}
  Given $n$ we can find $\Pi^0_1$ singletons $X,Y$ such that $X\nleT
  Y^{(n)}$ and $Y\nleT X^{(n)}$.
\end{thm}

\begin{proof}
  Let $X_n,Y_n$ be such that $0^{(n)}\leT X_n\leT0^{(n+1)}$ and
  $0^{(n)}\leT Y_n\leT0^{(n+1)}$ and $X_n\nleT Y_n$ and $Y_n\nleT
  X_n$.  Note that $X_n$ and $Y_n$ are $\Delta^{0,0^{(n)}}_2$ and
  hence $\Pi^{0,0^{(n)}}_2$ singletons.  Therefore, by the proof of
  Lemma \ref{lem:0} we may safely assume that $X_n$ and $Y_n$ are
  $\Pi^{0,0^{(n)}}_1$ singletons.  Apply Lemma \ref{lem:Pn} to
  $P_n=\{X_n,Y_n\}$ to get $X_0=H^n_0(X_n)$ and $Y_0=H^n_0(Y_n)$.
  Note that $P_0=\{X_0,Y_0\}$ is a $\Pi^0_1$ set, hence $X_0$ and
  $Y_0$ are $\Pi^0_1$ singletons.  Since $X_n\nleT
  Y_n\oplus0^{(n)}\eqT Y_0^{(n)}$ and $X_n\oplus0^{(n)}\eqT
  X_0\oplus0^{(n)}$ we have $X_0\nleT Y_0^{(n)}$, and similarly
  $Y_0\nleT X_0^{(n)}$.  Letting $X=X_0$ and $Y=Y_0$ we obtain our
  theorem.
\end{proof}

\begin{thm}
  Given $n$ we can find a countable $\Pi^0_1$ set $P$ such that some
  $Z\in P$ is not a $\Pi^0_n$ singleton.
\end{thm}

\begin{proof}
  Let $P_n$ be a countable $\Pi^0_1$ set such that some $Z_n\in P_n$
  is not isolated in $P_n$.  Treating $P_n$ as a $\Pi^{0,0^{(n)}}_1$
  set, apply Lemma \ref{lem:Pn} and note that $P_0$ is a countable
  $\Pi^0_1$ set and $Z_0=H^n_0(Z_n)$ is not isolated in $P_0$.  We
  claim that $Z_0$ is not a $\Pi^0_n$ singleton.  Otherwise, let $e$
  be such that $\{Z_0\}=\{X\mid e\notin X^{(n)}\}$.  Since $e\notin
  Z_0^{(n)}$ and $Z_0\in P_0$ and $X_0^{(n)}\eqT X_n\oplus0^{(n)}$
  uniformly for all $X_n\in P_n$ and $X_0=H^n_0(X_n)\in P_0$, there
  exists $j$ such that $e\notin X_0^{(n)}$ for all $X_n\in P_n$ such
  that $X_n\res j=Z_n\res j$.  But $Z_n$ is not isolated in $P_n$, so
  there exists $X_n\in P_n$ such that $X_n\res j=Z_n\res j$ and
  $X_n\ne Z_n$.  Thus $e\notin X_0^{(n)}$ and $X_0\ne Z_0$, a
  contradiction.  Letting $P=P_0$ and $Z=Z_0$ we obtain our theorem.
\end{proof}

\section{The full versions}
\label{sec:full}

In order to prove the full versions of Harrington's theorems, we need
to show that Lemma \ref{lem:Pn} holds with $n$ replaced by $\omega$.
To this end we first draw out some effective uniformities which are
implicit in the proofs of Lemmas \ref{lem:0} and \ref{lem:PQ}.

\begin{ntn}
  Let $W_e^A$ for $e=0,1,2,\ldots$ be a standard enumeration of all
  $A$-recursively enumerable subsets of $\NN^*$.  Then
  \begin{center}
    $T_e^A=\{\sigma\in\NN^*\mid(\forall n\le|\sigma|)\,(\sigma\res
    n\notin W_e^A)\}$
  \end{center}
  for $e=0,1,2,\ldots$ is a standard enumeration of all $\Pi^{0,A}_1$
  trees.  Hence $P_e^A=[T_e^A]$ for $e=0,1,2,\ldots$ is a standard
  enumeration of all $\Pi^{0,A}_1$ sets.
\end{ntn}

\begin{rem}
  \label{rem:unif}
  If $F$ is an $A$-recursive treemap and $T$ is a $\Pi^{0,A}_1$ tree,
  then $F(T)$ is again a $\Pi^{0,A}_1$ tree.  Moreover, this holds
  uniformly in the sense that there is a primitive recursive function
  $f$ such that $T_{f(e)}^A=F(T_e^A)$ and $P_{f(e)}^A=F(P_e^A)$ for
  all $e$, and we can compute a primitive recursive index of $f$
  knowing only an $A$-recursive index of $F$.
\end{rem}

The next two lemmas are refinements of Lemmas \ref{lem:0} and
\ref{lem:PQ} respectively.

\begin{lem}[refining Lemma \ref{lem:0}]
  \label{lem:0refined}
  There is a primitive recursive function $f$ with the following
  property.  Given $e$ we can effectively find an $A$-recursive
  treemap $F:T_e^{A'}\to T_{f(e)}^A$ which induces a homeomorphism
  $F:P_e^{A'}\cong P_{f(e)}^A$.  It follows that $X\oplus A\eqT
  F(X)\oplus A$ uniformly for all $X\in P_e^{A'}$.
\end{lem}

\begin{proof}
  Let $T=T_e^{A'}$ and $P=P_e^{A'}$.  Since $T_e^{A'}$ is uniformly
  $\Pi^{0,A'}_1$, it is uniformly $\Pi^{0,A}_2$, say
  $T=T_e^{A'}=\{\sigma\mid\forall i\,\exists j\,R(\sigma,e,i,A\res
  j)\}$ where $R\subseteq\NN^*\times\NN\times\NN\times\NN^*$ is a
  fixed primitive recursive predicate.  Let $(-,-)$ be a fixed
  primitive recursive one-to-one mapping of $\NN\times\NN$ onto $\NN$
  such that $m\le(m,n)$ and $n\le(m,n)$ for all $m$ and $n$.  Define
  $Q=[\widehat{T}]$ where
  $\widehat{T}=\{\sigma\oplus\tau\mid|\sigma|=|\tau|$ and
  $(\forall(n,i)<|\tau|)\,(\tau((n,i))=$ the least $j$ such that
  $R(\sigma\res n,e,i,A\res j))\}$.  Thus $Q=\{X\oplus\widehat{X}\mid
  X\in P\}$ where $\widehat{X}((n,i))=$ the least $j$ such that
  $R(X\res n,e,i,A\res j)$.  Moreover, we have an $A$-recursive
  treemap $F:T\to\widehat{T}$ given by
  $F(\sigma)=\sigma\oplus\widehat{\sigma}$ for all $\sigma\in T$,
  where $|\sigma|=|\widehat{\sigma}|$ and
  $(\forall(n,i)<|\sigma|)\,(\widehat{\sigma}((n,i))=$ the least $j$
  such that $R(\sigma\res n,e,i,A\res j))$.  Although we cannot expect
  to have $F(T)=\widehat{T}$, we nevertheless have
  $F:[T]\cong[\widehat{T}]$, i.e., $F:P\cong F(P)=Q$, and
  $F(X)=X\oplus\widehat{X}$ and $X\oplus A\eqT F(X)\oplus A$ uniformly
  for all $X\in P$.  The definition of $\widehat{T}$ shows that
  $\widehat{T}$ is uniformly $A$-recursive, hence uniformly
  $\Pi^{0,A}_1$, so we can find a fixed primitive recursive function
  $f$ such that $T_{f(e)}^A=\widehat{T_e^{A'}}$ for all $e$ and $A$.
\end{proof}

\begin{lem}[refining Lemma \ref{lem:PQ}]
  \label{lem:PQrefined}
  There is a primitive recursive function $h$ with the following
  property.  Given $e$ we can effectively find an $A'$-recursive
  treemap $H:T_e^{A'}\to T_{h(e)}^A$ which induces a homeomorphism
  $H:P_e^{A'}\cong P_{h(e)}^A$ such that $X\oplus A'\eqT H(X)\oplus
  A'\eqT(H(X)\oplus A)'$ uniformly for all $X\in P_e^{A'}$.
\end{lem}

\begin{proof}
  Let $G$ be the specific $A'$-recursive treemap which was constructed
  in the proof of Lemma \ref{lem:PQ}.  By Remark \ref{rem:unif} we can
  find a primitive recursive function $g$ such that for all $e$ we
  have $G(T_e^{A'})=T_{g(e)}^{A'}$ and the restriction of $G$ to
  $T_e^{A'}$ is a treemap from $T_e^{A'}$ to $T_{g(e)}^{A'}$ which
  induces a homeomorphism $G:P_e^{A'}\cong P_{g(e)}^{A'}$.  By
  construction of $G$ we have $X\oplus A'\eqT G(X)\oplus
  A'\eqT(G(X)\oplus A)'$ uniformly for all $X\in P_e^{A'}$.  Now
  applying Lemma \ref{lem:0refined} we obtain an $A$-recursive treemap
  $F:T_{g(e)}^{A'}\to T_{f(g(e))}^A$ which induces a homeomorphism
  $F:P_{g(e)}^{A'}\cong P_{f(g(e))}^A$ such that $Y\oplus A\eqT
  F(Y)\oplus A$ uniformly for all $Y\in P_{g(e)}^A$.  Thus the treemap
  $H=F\circ G:T_e^{A'}\to T_{f(g(e))}^A$ induces a homeomorphism
  $F\circ G=H:P_e^{A'}\cong P_{f(g(e))}^A$ such that $X\oplus A'\eqT
  H(X)\oplus A'\eqT(H(X)\oplus A)'$ uniformly for all $X\in P_e^{A'}$.
  Our lemma follows upon defining $h(e)=f(g(e))$.
\end{proof}

We now show that Lemma \ref{lem:Pn} holds with $n$ replaced by
$\omega$.

\begin{lem}
  \label{lem:Pomega}
  Given a $\Pi^{0,0^{(\omega)}}_1$ set $P_\omega$ we can effectively
  find a $\Pi^0_1$ set $P_0$ and a homeomorphism
  $H^\omega_0:P_\omega\cong P_0$ such that $X_\omega\oplus
  0^{(\omega)}\eqT X_0\oplus 0^{(\omega)}\eqT X_0^{(\omega)}$
  uniformly for all $X_\omega\in P_\omega$ and
  $X_0=H^\omega_0(X_\omega)\in P_0$.
\end{lem}

\begin{proof}
  Recall that $0^{(\omega)}=\{(i,n)\mid i\in0^{(n)}\}$.  Since
  $P_\omega$ is a $\Pi^{0,0^{(\omega)}}_1$ set, Remark \ref{rem:Pi01A}
  implies the existence of a tree $T_\omega\leT0^{(\omega)}$ such that
  $P_\omega=[T_\omega]$ and $\{\sigma\mid|\sigma|\le n$, $\sigma\in
  T_\omega\}\leT0^{(n)}$ uniformly for all $n$.  Define
  \begin{center}
    $T_{e,n}=\{\sigma\mid|\sigma|\le n\}\cup\{\sigma\mid|\sigma|>n$,
    $\sigma\res n\in T_\omega$, $\sigma\in T_e^{\langle
      n\rangle\cat0^{(n)}}\}$.
  \end{center}
  Thus $T_{e,n}$ is a $\Pi^{0,0^{(n)}}_1$ tree, hence
  $P_{e,n}=[T_{e,n}]$ is a $\Pi^{0,0^{(n)}}_1$ set, uniformly in $n$.

  In the vein of Lemma \ref{lem:PQrefined}, we claim there is a
  primitive recursive function $k$ with the following property.  Given
  $e$ and $n$ we can effectively find a $0^{(n+1)}$-recursive treemap
  \begin{center}
    $H_{e,n}:T_{e,n+1}\to T_{k(e),n}$
  \end{center}
  which induces a homeomorphism $H_{e,n}:P_{e,n+1}\cong P_{k(e),n}$
  such that $X\oplus0^{(n+1)}\eqT
  H_{e,n}(X)\oplus0^{(n+1)}\eqT(H_{e,n}(X)\oplus0^{(n)})'$ uniformly
  for all $X\in P_{e,n+1}$, and in addition $H_{e,n}(\sigma)=\sigma$
  for all $\sigma$ such that $|\sigma|\le n$.

  To prove our claim, let $r$ be a 3-place primitive recursive
  function such that
  $T_{r(e,n,\sigma)}^{0^{(n)}}=\{\tau\mid\sigma\cat\tau\in T_{e,n}\}$
  for all $e,n,\sigma$.  We can then write
  \begin{center}
    $T_{e,n+1}=\{\sigma\mid|\sigma|\le
    n\}\cup\{\sigma\cat\tau\mid|\sigma|=n$, $\tau\in
    T_{r(e,n+1,\sigma)}^{0^{(n+1)}}\}$.
  \end{center}
  Since $n$ is uniformly computable from $\langle
  n\rangle\cat0^{(n)}$, we can find a primitive recursive function $k$
  such that
  \begin{center}
    $T_{k(e),n}=\{\sigma\mid|\sigma|\le
    n\}\cup\{\sigma\cat\tau\mid|\sigma|=n$, $\tau\in
    T_{h(r(e,n+1,\sigma))}^{0^{(n)}}\}$
  \end{center}
  where $h$ is as in Lemma \ref{lem:PQrefined}.  For all $\sigma$ and
  $\tau$ such that $|\sigma|=n$ and $\tau\in
  T_{r(e,n+1,\sigma)}^{0^{(n+1)}}$ let
  $H_{e,n}(\sigma\cat\tau)=\sigma\cat H(\tau)$ where
  $H:T_{r(e,n+1,\sigma)}^{0^{(n+1)}}\to
  T_{h(r(e,n+1,\sigma))}^{0^{(n)}}$ is as in Lemma
  \ref{lem:PQrefined}.  Clearly $k(e)$ and $H_{e,n}$ have the required
  properties, so our claim is proved.

  Let $k$ and $H_{e,n}$ be as in the above claim.  By the Recursion
  Theorem (see Rogers \cite[Chapter 11]{rogers}) let $e$ be a
  \emph{fixed point} of $k$, so that $T_{k(e)}^A=T_e^A$ for all $A$,
  hence $T_{k(e),n}=T_{e,n}$ for all $n$.  Using this $e$ define
  $H_n=H_{e,n}$ and $T_n=T_{e,n}$ and $P_n=P_{e,n}=[T_n]$ for all $n$.
  As in the proof of Lemma \ref{lem:Pn} we have uniformly for each
  $s>n$ a $0^{(s)}$-recursive treemap $H^s_n=H_n\circ\cdots\circ
  H_{s-1}:T_s\to T_n$ which induces a homeomorphism $H^s_n:P_s\cong
  P_n$ such that $X\oplus0^{(s)}\eqT
  H^s_n(X)\oplus0^{(s)}\eqT(H^s_n(X))^{(s-n)}$ uniformly for all $X\in
  P_s$, and in addition $H^s_n(\sigma)=\sigma$ for all $\sigma$ such
  that $|\sigma|\le n$.  We also have for each $n$ a
  $0^{(\omega)}$-recursive treemap $H^\omega_n:T_\omega\to T_n$ which
  induces a homeomorphism $H^\omega_n:P_\omega\cong P_n$, namely
  $H^\omega_n(\sigma)=H^{|\sigma|}_n(\sigma)$ if $|\sigma|>n$ and
  $H^\omega_n(\sigma)=\sigma$ if $|\sigma|\le n$.  Note also that for
  all $n<s<t<\omega$ we have $H^t_n=H^s_n\circ H^t_s$ and
  $H^\omega_n=H^s_n\circ H^\omega_s$.  Finally, given $X_\omega\in
  P_\omega$ let $X_n=H^\omega_n(X_\omega)$ for all $n$.  Then
  $X_\omega\res n=X_n\res n$ and $X_n\oplus0^{(n)}\eqT
  X_0\oplus0^{(n)}\eqT X_0^{(n)}$ uniformly for all $n$ and all
  $X_\omega\in P_\omega$, hence $X_\omega\oplus0^{(\omega)}\eqT
  X_0\oplus0^{(\omega)}\eqT X_0^{(\omega)}$ uniformly for all
  $X_\omega\in P_\omega$.  This completes the proof.
\end{proof}

We now present Harrington's construction of arithmetically
incomparable arithmetical singletons.

\begin{thm}
  \label{thm:incomp}
  There is a pair of arithmetically incomparable $\Pi^0_1$ singletons.
\end{thm}

\begin{proof}
  Let $X_\omega,Y_\omega$ be such that $0^{(\omega)}\leT
  X_\omega\leT0^{(\omega+1)}$ and $0^{(\omega)}\leT
  Y_\omega\leT0^{(\omega+1)}$ and $X_\omega\nleT Y_\omega$ and
  $Y_\omega\nleT X_\omega$.  Note that $X_\omega$ and $Y_\omega$ are
  $\Delta^{0,0^{(\omega)}}_2$ and hence $\Pi^{0,0^{(\omega)}}_2$
  singletons.  Therefore, by the proof of Lemma \ref{lem:0} we may
  safely assume that $X_\omega$ and $Y_\omega$ are
  $\Pi^{0,0^{(\omega)}}_1$ singletons.  Apply Lemma \ref{lem:Pomega}
  to $P_\omega=\{X_\omega,Y_\omega\}$ to get a $\Pi^0_1$ set $P_0$ and
  a homeomorphism $H^\omega_0:P_\omega\cong P_0$.  Let
  $X_0=H^\omega_0(X_\omega)$ and $Y_0=H^\omega_0(Y_\omega)$.  Since
  $P_0=\{X_0,Y_0\}$ it follows that $X_0$ and $Y_0$ are $\Pi^0_1$
  singletons.  Since $X_\omega\nleT Y_\omega\oplus0^{(\omega)}\eqT
  Y_0^{(\omega)}$ and $X_\omega\oplus0^{(\omega)}\eqT
  X_0\oplus0^{(\omega)}$ we have $X_0\nleT Y_0^{(\omega)}$, and
  similarly $Y_0\nleT X_0^{(\omega)}$.  In particular $X_0$ and $Y_0$
  are arithmetically incomparable, Q.E.D.
\end{proof}

Finally we present Harrington's construction of a ranked point which
is not an arithmetical singleton.  This refutes a conjecture which had
been known as McLaughlin's Conjecture.  Note that McLaughlin's
Conjecture was natural in view of Remark \ref{rem:tanaka} above.

\begin{thm}
  \label{thm:mclaughlin}
  There is a countable $\Pi^0_1$ set $P$ such that some $Z\in P$ is
  not an arithmetical singleton.
\end{thm}

\begin{proof}
  Let $P_\omega$ be a countable $\Pi^0_1$ set such that some
  $Z_\omega\in P_\omega$ is not isolated in $P_\omega$.  Apply Lemma
  \ref{lem:Pomega} and note that $P_0$ is a countable $\Pi^0_1$ set
  and $Z_0=H^\omega_0(Z_\omega)\in P_0$ is not isolated in $P_0$.  We
  claim that $Z_0$ is not an arithmetical singleton.  Otherwise, let
  $e$ be such that $\{Z_0\}=\{X\mid e\in X^{(\omega)}\}$.  Since $e\in
  Z_0^{(\omega)}$ and $Z_0\in P_0$ and $X_0^{(\omega)}\eqT
  X_\omega\oplus0^{(\omega)}$ uniformly for all $X_\omega\in P_\omega$
  and $X_0=H^\omega_0(X_\omega)\in P_0$, there exists $j$ such that
  $e\in X_0^{(\omega)}$ for all $X_\omega\in P_\omega$ such that
  $Z_\omega\res j\subset X_\omega$.  But $Z_\omega$ is not isolated in
  $P_\omega$, so there exists $X_\omega\in P_\omega$ such that
  $Z_\omega\res j\subset X_\omega$ and $X_\omega\ne Z_\omega$.  Thus
  $e\in X_0^{(\omega)}$ and $X_0\ne Z_0$, a contradiction.  Letting
  $P=P_0$ and $Z=Z_0$ we obtain our theorem.
\end{proof}

\begin{rem}
  \label{rem:alpha}
  Modifying the proof of Lemma \ref{lem:Pomega}, it is easy to replace
  $\omega$ by a small recursive ordinal such as $\omega+\omega$ or
  $\omega\cdot\omega$ or $\omega^\omega$.  Harrington
  \cite{harrington-mclaughlin} and Gerdes \cite{gerdes-nusm} have
  shown that Lemma \ref{lem:Pomega} and consequently Theorems
  \ref{thm:incomp} and \ref{thm:mclaughlin} hold generally with
  $\omega$ replaced by any recursive ordinal.
\end{rem}

\bibliographystyle{plain}
\bibliography{misc}

\begin{thebibliography}{1}

\bibitem{masshyp}
Joshua~A. Cole and Stephen~G. Simpson.
\newblock Mass problems and hyperarithmeticity.
\newblock {\em Journal of Mathematical Logic}, 7:125--143, 2008.

\bibitem{feferman-generic}
Solomon Feferman.
\newblock Some applications of the notions of forcing and generic sets.
\newblock {\em Fundamenta Mathematicae}, 56:325--345, 1965.

\bibitem{gerdes-nusm}
Peter~M. Gerdes.
\newblock Harrington's solution to {M}c{L}aughlin's conjecture and non-uniform
  self-moduli.
\newblock Preprint, 27 pages, 15 December 2010, arXiv:1012.3427v1, submitted
  for publication.

\bibitem{harrington-arith}
Leo Harrington.
\newblock Arithmetically incomparable arithmetic singletons.
\newblock Handwritten, 28 pages, April 1975.

\bibitem{harrington-mclaughlin}
Leo Harrington.
\newblock Mc{L}aughlin's conjecture.
\newblock Handwritten, 11 pages, September 1976.

\bibitem{rogers}
Hartley Rogers{, Jr.}
\newblock {\em Theory of {R}ecursive {F}unctions and {E}ffective
  {C}omputability}.
\newblock McGraw-Hill, 1967.
\newblock XIX + 482 pages.

\bibitem{sacks-hrt}
Gerald~E. Sacks.
\newblock {\em Higher {R}ecursion {T}heory}.
\newblock Perspectives in {M}athematical {L}ogic. Springer-Verlag, 1990.
\newblock XV + 344 pages.

\bibitem{tanaka-arith}
Hisao Tanaka.
\newblock A property of arithmetic sets.
\newblock {\em Proceedings of the American Mathematical Society}, 31:521--524,
  1972.

\end{thebibliography}

\end{document}